\documentclass{eect}
\usepackage{amsmath}
  \usepackage{graphics} 
  \usepackage{epsfig} 
     \usepackage{graphicx,amssymb,color,cancel,lineno}
     \usepackage{empheq}
 \usepackage[colorlinks=true]{hyperref}
\hypersetup{urlcolor=blue, citecolor=red}

\def\currentvolume{2}
 \def\currentissue{4}
  \def\currentyear{2013}
   \def\currentmonth{December}
    \def\ppages{1--XX}
     \def\DOI{10.3934/eect.2013.2.xx}

\MHInternalSyntaxOn
\providecommand*\phantomword[3][c]{%
\mathchoice
{\MT_phantom_word:NNnn #1\displaystyle {#2}{#3}}%
{\MT_phantom_word:NNnn #1\textstyle {#2}{#3}}%
{\MT_phantom_word:NNnn #1\scriptstyle {#2}{#3}}%
{\MT_phantom_word:NNnn #1\scriptscriptstyle {#2}{#3}}%
}
\def\MT_phantom_word:NNnn #1#2#3#4{%
\@begin@tempboxa\hbox{$\m@th#2#4$}%
\setlength\@tempdima{\widthof{$\m@th#2#3$}}%
\hbox{\hb@xt@\@tempdima{\csname bm@#1\endcsname}}%
\@end@tempboxa}
\MHInternalSyntaxOff

\def\mc{\mathcal}

\newtheorem{theorem}{Theorem}[section]

\newtheorem{lemma}[theorem]{Lemma}

\theoremstyle{definition}

\title[Evolution Equations and Control Theory]
      {Uniform stabilization of a\\ multilayer Rao-Nakra sandwich beam}

\author[A. \"{O}zkan \"{O}zer and Scott W. Hansen]{}


 \email{aozer@uwaterloo.ca}
 \email{shansen@iastate.edu}

\thanks{Research of the second author was supported in part by the National Science Foundation under grant DMS- 1312952}

\begin{document}
\maketitle

\centerline{\scshape A. \"{O}zkan \"{O}zer}
\medskip
{\footnotesize
 \centerline{Department of Applied Mathematics, University of Waterloo}
   \centerline{Waterloo, ON N2L3G1, Canada}
} 

\medskip

\centerline{\scshape Scott W. Hansen}
\medskip
{\footnotesize
 \centerline{Department of Mathematics, Iowa State University}
   \centerline{Ames, Iowa 50011, USA}
}



\begin{abstract}
We consider the problem of boundary feedback stabilization of a multilayer Rao-Nakra sandwich beam.  We  show that the eigenfunctions of the decoupled system form a Riesz basis. This allows us to deduce that the decoupled system  is exponentially stable. Since the coupling terms are compact, the exponential stability of the coupled system  follows from the strong stability of the coupled system, which is proved using a unique continuation
result for  the overdetermined homogenous system in the case of zero feedback.
\end{abstract}

\section{Introduction}
A {\it sandwich beam}  is an engineering model for a  three-layer  beam consisting of stiff  outer face plates and a more compliant   inner core layer.
Sandwich beam models found in the literature include the models of Mead and Markus \cite{Mead-Marcus},  Rao and Nakra  (RN)  \cite{Rao-Nakra},
Yan and Dowell \cite{Yan-Dowell} and others.    The RN model assumes continuous, piecewise linear displacements through the cross-sections, with the Kirchhoff hypothesis imposed on the face plates.   Transverse, longitudinal and rotational inertial forces are included in the modeling.
   In \cite{Hansen}  several possible multilayer generalizations of the basic three layer sandwich beam structure are derived and analyzed (in the form of multilayer plates).     In this paper we consider a multilayer generalization of the RN model described in \cite{Hansen}.
The model consists of  $2m+1$ alternating stiff and complaint (core) layers, with stiff layers on  the outside.   The stiff layers assume the Kirchhoff hypothesis, while the compliant layers admit shear.      The equations of motion for the associated beam model can be written:
\begin{equation}\left\{ \begin{array}{l}
\ddot z -\alpha \ddot z'' +  K  z'''' -   N^T {\bf{h}}_E {\bf{ G}}_E \phi_E' = 0 ~~~ {\rm{on}}~~ \Omega\times \mathbb{R}^+ \\
   {\bf{h}}_\mc{O} {\bf{p}}_\mc{O} {\ddot v}_\mc{O} -{\bf{h}}_\mc{O} {\bf{E}}_\mc{O}  {v}_\mc{O}'' + {\bf{B}}^T  {\bf{ G}}_E \phi_E  = 0~~ {\rm{on}} ~~ \Omega\times \mathbb{R}^+ \\
 {\rm where}~( {\bf{B}} { v}_\mc{O}={\bf{h}}_E \phi_E-{\bf{h}}_E  N z')
 \end{array} \right.
\label{main}
\end{equation}
where $\Omega=(0,L),$ primes denote differentiation with respect to the spatial variable $x$
and dots denote differentiation with respect to time $t.$

In the above, $z$ represents the transverse displacement, $\phi^i$ denotes the shear angle in the $i^{\rm{th}}$ layer, $\phi_E=[\phi^2,\phi^4,\ldots,\phi^{2m} ]^{\rm T},$ $v^i$ denote the longitudinal displacement
along the center of the   $i^{\rm{th}}$   layer, and $v_{\mc O}=[v^1,v^3, \ldots, v^{2m+1}]^{\rm T}.$      Throughout this paper we use the convention that
quantities relating to the  stiff layers have odd indices $1,3,\ldots 2m+1$ and quantities relating to the even layers have even indices $2,4,\ldots 2m$.
    In addition, $m,\alpha, K$ are positive physical constants, and
\begin{eqnarray}
\nonumber &{\bf{p}}_{\mc O}={\rm{diag}}~ (\rho_1,  \ldots,\rho_{2m+1}), ~~{\bf{h}}_E={\rm{diag}}~(h_2, \ldots, h_{2m}),  ~~{\bf{h}}_{\mc O}={\rm{diag}}~(h_1, \ldots, h_{2m+1})&\\
\nonumber &{\bf{ G}}_E={\rm{diag}}~( G_2, \ldots, G_{2m}), ~~{\bf{E}}_{\mc O}={\rm{diag}}~(E_1,  \ldots, E_{2m+1})&
\end{eqnarray}
where $\rho_i, h_i, E_i,  G_i$ denote the density, thickness, Young's modulus, and shear modulus of the $i^{\rm{th}}$ layer, respectively. The vector $ N$ is defined as
$N={\bf{h}}_E^{-1}{\bf{A }}{\bf{h}}_{\mc O}  \vec 1_{\mc O} +  \vec 1_E$
where ${\bf{A}}=(a_{ij})$  and ${\bf{B}}=(b_{ij})$ are the $m\times(m+1)$ matrices
$$a_{ij}  = \left\{ \begin{array}{l}
1/2,~~{\rm{  if   }}~~j = i~~{\rm{  or }}~~j = i + 1 \\
~~0,\quad{\rm{   otherwise}} \\
\end{array} \right., ~~b_{ij}=\left\{ \begin{array}{l}
(-1)^{i+j+1},~~{\rm{  if   }}~~j = i~~{\rm{  or }}~~j = i + 1 \\
~~0, \quad\quad\quad\quad {\rm{    otherwise}} \\
\end{array} \right. $$
and $ \vec 1_{\mc O}$ and $\vec 1_E$ denote column vectors with all entries of 1 in $\mathrm{R}^{m+1}$ and $\mathrm{R}^{m},$ respectively.

The aim of this paper is to prove that the uniform exponential stability of the RN system with standard boundary damping applied at one end point.
Consider (\ref{main}) with the following boundary conditions
\begin{eqnarray} \left\{ \begin{array}{l}  z(0,t)=z'(0,t)= z(L,t)=0,~ z''(L,t) + \gamma_0 ~\dot z'(L,t)=0\label{bdry1}\\
 {v}_\mc{O}(0,t) = 0,  ~ {v}_\mc{O}'(L,t) + {{\Upsilon}_{\mc \mc{O}}}~\dot{v}_\mc{O}(L,t)=0,  \quad\quad\quad\quad\quad\quad\quad\quad  {\rm{on}} \quad\mathbb{R}^+ \label{bdry2}
  \end{array} \right.
\end{eqnarray}
and the initial conditions
\begin{eqnarray}
 z(x,0)=z^0,  ~\dot z(x,0)=z^1, ~ {v}_{\mc O}(x,0)= {v}^0_{\mc{O}}, ~ {\dot v}_{\mc{O}}(x,0)= {v}^1_\mc{O}, \quad  {\rm{on}} \quad \Omega\label{initial}
\end{eqnarray}
where  ${\Upsilon}_\mc{O}= {\rm diag}
~(\gamma_1, \gamma_3 \cdots, \gamma_{2m+1}),$ and  $\gamma_i \in \mathbb{R}^+, i=0,1,3,\ldots,2m+1$ denote constant positive feedback gains. Throughout the paper, we assume
\begin{eqnarray}\label{assump}~\sqrt{\frac{\alpha}{K}}\ne \gamma_0\quad {\rm{and}}\quad \sqrt{\frac{\rho_k}{E_k}}\ne \gamma_k \quad {\rm{for}} \quad k=1,3,\ldots, 2m+1.\end{eqnarray}

\subsection{Background}
 Boundary controllability of   (\ref{main}) has been studied in several papers.  For the three layer case,   in \cite{Rajaram-Hansen2}  the multiplier method  was used to prove exact controllability with a control for each equation applied at  an end point.    The moment method was used in \cite{Rajaram-Hansen3} to obtain boundary controllability for the multilayer case, but with the condition that wave speeds of the layers be distinct.    The same approach was used to prove simultaneous controllability (i.e., with one boundary control instead of three) for the three layer case in \cite{Rajaram-Hansen1}.     In \cite{O-Hansen3} exact boundary controllability of  the general multilayer system was proved for a variety of boundary conditions: clamped, hinged, clamped-hinged, and hinged-clamped. The results in \cite{O-Hansen3}  improve earlier results in that
 there are no restrictions on the wave speeds or  the size of ${{\bf{ G}}}$ and moreover, exact controllability is proved in the optimal time (determined by characteristics).    In \cite{Hansen-Oleg}, \cite{Hansen-Oleg1}  exact controllability results for the multilayer RN plate system analogous to (\ref{main}) with locally distributed control in a neighborhood of a portion of the boundary were obtained by the method of Carleman estimates.

Stability results for layered beam systems closely related to (\ref{main}) subject to internal damping proportional to rate-of-shear  in one or more layers have been studied in several papers; \cite{AH1}, \cite{AH2}, \cite{H-L},    \cite{Wang-Guo}.    In particular,  the approach used in \cite{H-L} was successfully applied  in the dissertation \cite{AAA}  to obtain uniform exponential stability results the system (\ref{main}) with rate-of shear damping included in the
compliant layers.

Concerning boundary feedback stabilization of layered beam models,    spectral methods (based on the  Riesz basis property) are applied   in
\cite{Wang-Xu-Yung} to prove exponential stability results   for  a laminated beam model in \cite{H-S}.    A similar approach is used in  \cite{Guo-Chentouf}  for the Mead-Markus model described in \cite{Fabiano-Hansen}.  There are also several  results  concerning the boundary feedback stabilization of  a single Rayleigh beam equation e.g., \cite{B. Rao},   where  a uniform exponential decay result is  obtained for a
Rayleigh beam by means of  a compact perturbation argument and \cite{Guo}   where  the Riesz basis approach is used  to obtain a similar stabilization result.      Some related uniform stabilization results for the Kirchhoff plate  are proved in \cite{LT1}.



Our main result is the following.


\begin{theorem}\label{finaltheorem} Assume (\ref{assump}). Then the semigroup generated by $\mc{A}$ is exponentially stable in $\mc{H}$, i.e., $\exists M>0,$ $\mu<0$ such that
$\mc{E}(t)\le M e^{\mu t}\mc{E}(0).$
Moreover, $\mu=\sup\{~\mbox{Re}\,\lambda~|~ \lambda \in \sigma(\mc{A})\}.$
\end{theorem}

In the above, $\mc{A},$ $\varepsilon(t),$ and $\mc{H}$  are defined in (\ref{semigroupfor1}), (\ref{energy}), and (\ref{semigroupdom}) respectively.


   Our methodology in this paper is a combination of techniques used in \cite{B. Rao}, \cite{Trigg}, and \cite{Wang-Xu-Yung}.    The decoupled system
(i.e., (\ref{main}),  with  ${\bf{G}}_E \equiv 0$)     consists of a Rayleigh beam equation and  $(m+1)$ wave equations.    We prove that the decoupled system
has a Riesz basis of eigenfunctions and  obtain explicit asymptotic estimates  on the eigenvalues.    In particular,  the eigenvalues  of the decoupled system
asymptotically lie along a finite number of vertical lines in the left half plane.   We are able to prove that the family of eigenfunctions and generalized eigenfunctions of the decoupled system form a Riesz basis and consequently (see \cite{Trigg0}), the  {\it spectrum determined growth condition} holds.
This allows us to prove the exponential stability of the decoupled closed-loop system (see Theorem \ref{exp-ray}).   We mention that the exponential stabilization  for the portion of the uncoupled system corresponding to the wave equations  is well-known results, e.g., \cite{Chen}, \cite{K-Z}, \cite{Lagnese}, \cite{LT2}.   Furthermore  a number of results are known for stabilization of the Rayleigh beam, e.g. \cite{Guo}, \cite{LT1}, \cite{B. Rao},    however none of these results  are
applicable  to the  Rayleigh beam with clamped-hinged boundary conditions which we consider.   Therefore we include a detailed proof of the exponential stability for the Rayleigh beam.
Next, we prove that the system (\ref{main})-(\ref{initial}) has a compact resolvent and  is a compact perturbation of the decoupled system.  Therefore, exponential stability of (\ref{main})-(\ref{initial}) follows from a perturbation theorem due to Triggiani \cite{Trigg} once it is shown that the semigroup generated by $\mc A$ is strongly stable (see Theorem \ref{stronglystable}).   Proving the strong stability involves use of dissipativity of the semigroup together with a nontrivial unique continuation argument that is proved
in \cite{O-Hansen3} in application to the associated boundary control problem.

    Our paper is organized as the following. In Section \ref{WP}, we give a semigroup formulation of (\ref{main})-(\ref{initial}). We prove that the semigroup is a $C_0-$semigroup of contractions on an appropriate Hilbert space. In Section \ref{Decoupled}, we first prove that the generalized functions corresponding to the single Rayleigh beam equation with the feedback applied to the moment forms a Riesz basis. Then, we show that the decoupled system, i.e. ${\bf{G}}_E \equiv 0$ in (\ref{main}), has the Riesz basis property. Finally, we show that the decoupled system is exponentially stable. In Section \ref{Coupled-Stab}, we prove that the system (\ref{main}) is a compact perturbation of the decoupled system. Finally, we prove our main stabilization result in Theorem \ref{finaltheorem}.

 \section{Well-posedness of the system}
 \label{WP}
Let
\begin{eqnarray}\nonumber U=:(u,{\bf{u}})^{\rm T}=(z, {v}_\mc{O})^{\rm T}, ~~~~V:= (v, {\bf{v}})^{\rm T}=(\dot z,  {\dot v}_\mc{O})^{\rm T}, ~~{\rm{and}} ~~ Y:=(U,V)^{\rm T}.
\end{eqnarray}
 Let also $L\varphi= \varphi - \alpha \varphi''$. From the Lax-Milgram theorem, $L:H^1_0(\Omega) \to H^{-1}(\Omega)$ is an isomorphism.  Then (\ref{main})-(\ref{initial}) can be formulated as\newpage
\begin{eqnarray}\nonumber &&\frac{dY}{dt}=\mc{A} Y:=\left( {\begin{array}{*{20}c}
     {0} & {I }  \\
   {A_1 } & {0}    \\
\end{array}} \right) \left( \begin{array}{c}
 U \\
 V \\
 \end{array} \right), ~~~~ Y(0)=(U(0),V(0))^{\rm T}=(z^0,  {v}^0_\mc{O}, z^1,  {v}^1_\mc{O} )^{\rm T}
\label{semigroupfor1}\end{eqnarray}
where
\begin{eqnarray}\nonumber A_1 U:=\left( {\begin{array}{*{20}c}
   L^{-1}\left( -Ku''''+ N^{\rm T} {\bf{h}}_E {\bf{G}}_E ~({\bf{h}}_E^{-1} {\bf{B}} {\bf{u}}'+  N u'') \right)  \\
   { {\bf{h}}^{-1}_{\mc O} {\bf{p}}^{-1}_{\mc O} } \left({ {\bf{h}}_{\mc O} {\bf{E}}_{\mc O} }\textbf{u}''-{\bf{B}}^{\rm T} {\bf{G}}_E~( {\bf{h}}^{-1}_E {\bf{B}} {\bf{u}}+ N u')\right) \\
\end{array}} \right).
\end{eqnarray}
Let $\left<u,v\right>_{\Omega}=\int_{\Omega} u \cdot \overline{v} ~dx$ where $u$ and $v$ may be scalar or vector valued. Define the bilinear forms $a$ and $c$ by
\begin{eqnarray}
\nonumber  c(z, v_\mc{O}; \hat{z},  {\hat{v}_\mc{O}})&=& \left<z,\hat{z}\right>_{\Omega}+\alpha\left<z',\hat{z}'\right>_{\Omega}+ \left< {\bf{h}}_\mc{O} {\bf{p}}_\mc{O} { v}_\mc{O},   {\hat{ v}}_\mc{O}\right>_{\Omega}\\
  \nonumber a(z, v_\mc{O}; \hat{z},  {\hat{v}_\mc{O}})&=& K  \left<z'',\hat{z}''\right>_{\Omega}+ \left< {\bf{h}}_\mc{O} {\bf{E}}_\mc{O} { v}_\mc{O}',   { \hat{v}}_\mc{O}'\right>_{\Omega}+   \left<{\bf{G}}_E {\bf{h}}_E \phi_E,{\hat{\phi}}_E\right>_{\Omega}\\
\nonumber &=& K  \left<z'',\hat{z}''\right>_{\Omega}+ \left< {\bf{h}}_\mc{O} {\bf{E}}_\mc{O} { v}_\mc{O}',   { \hat{v}}_\mc{O}'\right>_{\Omega}\\
 \label{sesqui}&& ~+\left<{\bf{G}}_E {\bf{h}}_E^{-1}\left( {\bf{B}} { v}_{\mc O}+N z'\right),\left( {\bf{B}} {\hat{v}}_{\mc O}+N \hat z'\right)\right>_{\Omega}.
\end{eqnarray}
The natural energy of the beam is given by
 \begin{eqnarray}
\label{energy}
\mc{E}(t)= \frac{1}{2}\left( a(z, v_\mc{O})+c(\dot z, \dot v_\mc{O})\right),\quad {\rm for ~all}~ t\in \mathbb{R}^+
 \end{eqnarray}
where  $a(\cdot), c(\cdot)$ are the quadratic forms that agree with $a(\cdot~;~ \cdot), c(\cdot~;~ \cdot)$ on the diagonal.
We define the Hilbert space $\mc H$ by
 \begin{eqnarray}
\label{semigroupdom}
&\mc{H}=     {X} \times {Y}; &\\
 & X= H^2_\#(\Omega)\times \left(H^1_*(\Omega)\right)^{(m+1)},\quad Y=      H^1_0(\Omega)\times (L^2(\Omega))^{(m+1)}&\nonumber
\end{eqnarray}
with the  energy inner product
 \begin{eqnarray}\label{inner}
\left<Y,\hat Y\right>_{\mc{H}}= a(U;\widehat{U})+c(V;\widehat{V})
  \end{eqnarray}
  where
\begin{eqnarray}
\nonumber H^2_\#(\Omega)&=&\{ u \in \left(H^2(\Omega) \cap H_{0}^1(\Omega)\right)~:~ u'(0)=0\},\\
\nonumber H^1_*(\Omega)&=&\{ u \in H^1(\Omega) ~:~ u(0)=0\}.
 \end{eqnarray}
\vspace{0.1in}

\noindent
\textbf{Characterization of the domain  of $\mc A:$}  Let $(u_1, {\bf u}_1, v_1, {\bf v}_1)\in  \mc  H $ such that  $A_1 \left( {\begin{array}{*{20}c}
u_1 \\
{\bf u}_1 \\
\end{array}} \right) \in  Y= H^1_0(\Omega)\times (L^2(\Omega))^{(m+1)},$ $V\in X$ and assume that boundary conditions (\ref{bdry1})  hold.
Also  let
$ (f, {\bf f})^{\rm T}\in X= H^2_\#(\Omega)\times (H^1_*(\Omega))^{(m+1)}$.
A calculation shows that
\begin{align}\nonumber c\left(A_1 \left( {\begin{array}{*{20}c}
u_1 \\
{\bf u}_1 \\
\end{array}} \right) ; \left( {\begin{array}{*{20}c}
f \\
{\bf f} \\
\end{array}} \right)\right)=&c\left(\left( {\begin{array}{*{20}c}
   L^{-1}\left( -K u_1''''+ N^{\rm T} {\bf{h}}_E {\bf{G}}_E \Phi_E' \right)  \\
   { {\bf{h}}^{-1}_\mc{O} {\bf{p}}^{-1}_\mc{O} } \left({ {\bf{h}}_\mc{O} {\bf{E}}_\mc{O} }{\bf u}_1''-{\bf{B}}^{\rm T} {\bf{G}}_E\Phi_E\right) \\
\end{array}} \right); \left( {\begin{array}{*{20}c}
f \\
{\bf f} \\
\end{array}} \right)\right)\\
\nonumber =&\left<\left( {\begin{array}{*{20}c}
   -K u_1''''+ N^{\rm T} {\bf{h}}_E {\bf{G}}_E \Phi_E'  \\
{ {\bf{h}}_\mc{O} {\bf{E}}_\mc{O} }{\bf u}_1''-{\bf{B}}^{\rm T} {\bf{G}}_E\Phi_E\\
\end{array}} \right), \left( {\begin{array}{*{20}c}
f \\
{\bf f} \\
\end{array}} \right)\right>_{\Omega}\\
\nonumber =&   -  \left<{ {\bf{h}}_\mc{O} {\bf{E}}_\mc{O} }{\bf u}_1', {\bf f}'\right>_{\Omega}- \left< {\bf{G}}_E \Phi_E, {\bf{h}}_E  N f' + {\bf{B}} {\bf f}\right>_{\Omega} \quad\quad\\
\nonumber & -K \left<u_1'', f''\right>_{\Omega} -K \gamma_0 v_1'(L)f'(L) - { {\bf{h}}_\mc{O} {\bf{E}}_\mc{O} }{{\Upsilon}}_\mc{O} {\bf{v}}_1 \cdot {\bf{f}}(L).   \end{align}
In the above,  ${\bf{B}} {\bf u}_1={\bf{h}}_E \Phi_E-{\bf{h}}_E  N {u}_1'.$
 Hence,  using the definition of $a$ in (\ref{semigroupfor1}), we find that the following identity holds:
 \begin{eqnarray}\nonumber &&c\left(A_1 \left( {\begin{array}{*{20}c}
u_1 \\
{\bf u}_1 \\
\end{array}} \right) ; \left( {\begin{array}{*{20}c}
f \\
{\bf f} \\
\end{array}} \right)\right)+a\left(\left( {\begin{array}{*{20}c}
u_1 \\
{\bf u}_1 \\
\end{array}} \right); \left( {\begin{array}{*{20}c}
f \\
{\bf f} \\
\end{array}} \right)\right)\\
\label{eq21} && ~~~~~~~ =-K \gamma_0 v_1'(L)f'(L)-{ {\bf{h}}_\mc{O} {\bf{E}}_\mc{O} }{{\Upsilon}}_\mc{O}{\bf v}_1(L) \cdot {\bf{f}}(L)  \quad\forall\, (f, {\bf f})^{\rm T}\in X.
\end{eqnarray}
We use  (\ref{eq21}) as a basis for the  variational definition of  the space  $\mc{D}({\mc{A}}).$  More precisely,
\begin{eqnarray}
\nonumber \mc D(\mc A)= \{(U,V)\in X\times X: A_1 U\in Y\mbox{ and (\ref{eq21}) holds}\}.
\end{eqnarray}
 Note that the operator   $\mc{A}: \mc{D}(\mc{A})\subset \mc{H} \to \mc{H}$  is densely defined. 

\begin{lemma} \label{surjective}$I-\mc{A}: \mc{D}(\mc{A})\to \mc{H}$ is surjective, i.e. ${\rm Range}(I-\mc{A})=\mc{H}.$
\end{lemma}

\noindent\emph{Proof.} Let $Y_1=(u_1,{\bf{u}_1}, v_1, {\bf{v}_1})^{\rm T}.$  For given $Y_2=(u_2,{\bf{u}_2}, v_2, {\bf{v}_2})^{\rm T} \in  \mc{H} $ we want to prove the solvability of the system $(I-\mc{A})Y_1=Y_2$  for $Y_1\in \mc{D}({\mc{A}}).$ This is equivalent to prove the solvability of the following system in $\mc{D}({\mc{A}}):$
\begin{eqnarray}
\label{oz1} L^{-1}\left(- K u_1''''+ N^{\rm T} {\bf{h}}_E {\bf{G}}_E ({\bf{h}}_E^{-1} {\bf{B}} {\bf{u}_1}'+  N u_1'')\right) =&& v_1-v_2 \\
\label{oz2}  {\bf{h}}_\mc{O}^{-1}{\bf{p}}_\mc{O}^{-1}\left({ {\bf{h}}_\mc{O} {\bf{E}}_\mc{O} }\textbf{u}_1''-{\bf{B}}^{\rm T} {\bf{G}}_E( {\bf{h}}^{-1}_E {\bf{B}} {\bf{u}_1}+ N u_1')\right)= &&{\bf{v}_1}-{\bf{v}_2}\\
\label{rc1}   u_1-v_1 =&& u_2 \\
\label{rangecon} {\bf{u}_1} - {\bf{v}_1}=&&   {\bf{u}_2}.
\end{eqnarray}
Let $(f,{\bf f})^{\rm T}\in H^2_\#(\Omega)\times (H^1_*(\Omega))^{(m+1)}. $ If we multiply (\ref{oz1}) by $f$ and dot product (\ref{oz2}) by ${\bf f}$ then  integrate by parts and  apply (\ref{rc1}), (\ref{rangecon})  we obtain
\begin{eqnarray}
 \label{rangecon2} \nonumber &&a\left( \left( {\begin{array}{*{20}c}
u_1 \\
{\bf u}_1 \\
\end{array}} \right) ;   \left( {\begin{array}{*{20}c}
f \\
{\bf f} \\
\end{array}} \right)\right) + c\left(\left( {\begin{array}{*{20}c}
u_1 \\
{\bf u}_1 \\
\end{array}} \right); \left( {\begin{array}{*{20}c}
f \\
{\bf f} \\
\end{array}} \right)\right)\\
\nonumber &&\quad\quad+ K \gamma_0 u_1'(L)f'(L)+{ {\bf{h}}_\mc{O} {\bf{E}}_\mc{O} }{{\Upsilon}}_\mc{O} {\bf u}_1(L) \cdot {\bf{f}}(L)\\
 && ~~= c\left(\left( {\begin{array}{*{20}c}
 u_2+v_2     \\
    {\bf{u}_2 + {\bf{v}}_2}\\
\end{array}} \right) ; \left( {\begin{array}{*{20}c}
f \\
{\bf f} \\
\end{array}} \right)\right) \nonumber\\ \label{rangecon3}  &&+ K \gamma_0 u_2'(L)f'(L) + { {\bf{h}}_\mc{O} {\bf{E}}_\mc{O} }{{\Upsilon}}_\mc{O}{\bf u}_2(L) \cdot {\bf{f}}(L)
\quad \forall   \left( {\begin{array}{*{20}c}
f \\
{\bf f} \\
\end{array}} \right) \in X. \end{eqnarray}
  The bilinear forms $a$ and $c$ are symmetric, bounded and coercive on $H^2_\#(\Omega)\times (H^1_*(\Omega))^{(m+1)}$ and  $H^1_0(\Omega)\times (L^2(\Omega))^{(m+1)},$ respectively. Moreover, the right hand side of (\ref{rangecon3}) is a bounded linear form on $H^2_\#(\Omega)\times (H^1_*(\Omega))^{(m+1)}$.    Therefore, by Lax-Milgram theorem, there exists a unique pair $(u_1, {\bf u}_1)^{\rm T}\in H^2_\#(\Omega)\times (H^1_*(\Omega))^{(m+1)}$ satisfying (\ref{rangecon3}). Therefore, (\ref{rangecon}) uniquely determines pair $(v_1, {\bf v}_1)^{\rm T}\in H^2_\#(\Omega)\times (H^1_*(\Omega))^{(m+1)}.$\\

The last step of our proof is to show that $Y_1\in \mc D(\mc A),$ i.e. $(u_1, {\bf u}_1, v_1, {\bf v}_1)^{\rm T}$ satisfies (\ref{eq21}). Now assume that $(f,{\bf f})^{\rm T}=( g, {{ \bf g}})^{\rm T}$ with $( g, {{ \bf g}})^{\rm T}\in C_0^{\infty}(\Omega).$ Then it follows that
\begin{eqnarray}
  \nonumber &&c\left( \left( {\begin{array}{*{20}c}
u_1 \\
{\bf u}_1 \\
\end{array}} \right) , \left( {\begin{array}{*{20}c}
g \\
{\bf g} \\
\end{array}} \right)\right) + a\left(\left( {\begin{array}{*{20}c}
u_1 \\
{\bf u}_1 \\
\end{array}} \right),\left( {\begin{array}{*{20}c}
g \\
{\bf g} \\
\end{array}} \right)\right) =c\left(\left( {\begin{array}{*{20}c}
 u_2+v_2     \\
    {\bf{u}_2 + {\bf{v}}_2}\\
\end{array}} \right) , \left( {\begin{array}{*{20}c}
g \\
{\bf g} \\
\end{array}} \right)\right)~~~~~ \end{eqnarray}
holds for all $( g, {{ \bf g}})^{\rm T}\in C_0^{\infty}(\Omega).$ Therefore in $(C_0^{\infty}(\Omega))'$ we have
\begin{eqnarray}\label{rangecon5}
\nonumber  L^{-1}\left(K u_1''''- N^{\rm T} {\bf{h}}_E {\bf{G}}_E ({\bf{h}}_E^{-1} {\bf{B}} {\bf{u}_1}'+  N u_1'')\right)  &&= u_1-u_2-v_2 \in H^1_0(\Omega)\\
  {\bf{h}}_\mc{O}^{-1}{\bf{p}}_\mc{O}^{-1} \left( -{ {\bf{h}}_\mc{O} {\bf{E}}_\mc{O} }\textbf{u}_1''+{\bf{B}}^{\rm T} {\bf{G}}_E( {\bf{h}}^{-1}_E {\bf{B}} {\bf{u}_1}+ N u_1')\right)&&= {\bf u_1}- {\bf{u}_2}-{\bf{v}_2}\in L^2(\Omega).\quad\quad
\end{eqnarray}
If we substitute (\ref{rangecon5}) in (\ref{rangecon3}) by setting  $v_1=u_1-u_2$ and ${\bf v}_1={\bf u}_1-{\bf u}_2$ we obtain that $Y_1$ satisfies (\ref{eq21}). This together with (\ref{rangecon5}) implies that $Y_1\in \mc{D}(\mc{A})$. $\hfill\square$

\begin{lemma} \label{disspppp}The infinitesimal generator $\mc{A}$ is dissipative on $\mc{H},$ and it satisfies  \begin{eqnarray}{\rm Re} \left<\mc{A}Y, Y\right>_{\mc{H}} =  -K \gamma_0 |v_1'(L)|^2 - {\bf{h}}_\mc{O} {\bf{E}}_\mc{O} {{\Upsilon}}_\mc{O}{\bf{v}}_1(L) \cdot{\bf{{\bar v}}}_1 (L)\le 0~~~~
\label{Dissp}
\end{eqnarray}
for all $Y=(u, {\bf{u}},v,{\bf{v}})^{\rm T}\in \mc D(\mc A).$
\end{lemma}

\noindent\emph{Proof.} By an easy calculation  we have the following
    \begin{align}
\nonumber \left<\mc{A}Y, Y\right>_{\mc{H}} =&   \left\{-K \left<u'', v''\right>_{\Omega} +K\left< v'', u''\right>_{\Omega}\right\}+\left\{ \left<{\bf{h}}_\mc{O} {\bf{E}}_\mc{O} {\bf{v}}', {\bf{u}'}\right>_{\Omega}-\left<{ {\bf{h}}_\mc{O} {\bf{E}}_\mc{O} }\textbf{u}', {\bf{v}}'\right>_{\Omega}\right\}~~~~\\
\nonumber   &+\left\{-\left< {\bf{G}}_E\left( {\bf{B}} {\bf{u}}' + {\bf{h}}_E   N u''\right) , {\bf{h}}^{-1}_E \left({\bf{B}} {\bf{v}'} + {\bf{h}}_E^{-1}N v''\right)\right>_{\Omega} \right. \\
 \nonumber &\left.\quad\quad\quad+ \left<{\bf{G}}_E  ({\bf{B}} {\bf{v}}'+  {\bf{h}}_E N v''), {\bf{h}}_E^{-1} ({\bf{B}} {\bf{u}}'+ {\bf{h}}_E N u'')\right>_{\Omega}\right\}\\
  \nonumber &-K \gamma_0|v_1'(L)|^2 - {\bf{h}}_\mc{O} {\bf{E}}_\mc{O} {{\Upsilon}}_\mc{O}{\bf{v}}_1(L) \cdot{\bf{{\bar v}}}_1 (L)\\
  \nonumber =&  -2i ~{\rm Im} \left\{~K \left<u'', v''\right>_{\Omega} +\left<{ {\bf{h}}_\mc{O} {\bf{E}}_\mc{O} }\textbf{u}', {\bf{v}}'\right>_{\Omega}\right\}\\
\nonumber   & -2i ~{\rm Im} \left< {\bf{G}}_E\left( {\bf{B}} {\bf{u}}' + {\bf{h}}_E   N u''\right) , {\bf{h}}^{-1}_E \left({\bf{B}} {\bf{v}'} + {\bf{h}}_E^{-1}N v''\right)\right>_{\Omega}  \\
  \nonumber &-K \gamma_0|v_1'(L)|^2 - {\bf{h}}_\mc{O} {\bf{E}}_\mc{O} {{\Upsilon}}_\mc{O}{\bf{v}}_1(L) \cdot{\bf{{\bar v}}}_1 (L).
\end{align}
Therefore  (\ref{Dissp}) follows. $\hfill\square$


\begin{lemma} \label{egigen0} The point spectrum of $\mc{A}$ does not contain $\lambda =0,$ i.e. the following eigenvalue problem
\begin{equation}\left\{ \begin{array}{l}
 K u'''' - N^{\rm T} {\bf{h}}_E {\bf{G}}_E \phi_E'=0 \\
  {\bf{h}}_\mc{O} {\bf{E}}_\mc{O}  {v}_\mc{O}'' - {\bf{B}}^T  {\bf{ G}}_E \phi_E=0 \\
   {\bf{B}} { v}_\mc{O}={\bf{h}}_E \phi_E-{\bf{h}}_E  N u'
 \end{array} \right.
\label{eq}
\end{equation}
with the boundary conditions
\begin{eqnarray}
 \label{bc} &&u(0,t)=u'(0,t)= u(L,t)=u''(L,t)=0, ~~v_{\mc O}(0)=v_{\mc O}'(L)=0
\end{eqnarray}
has only the trivial solution.
\end{lemma}


\noindent\emph{Proof.}  Let ${\bf T}=-{\rm{diag}}~(D_x^2,\ldots, D_x^2)$ be defined on the domain
${\rm Dom}({\bf T})=\{\psi\in (H^2(\Omega))^{m}~:~ \psi(0)=\psi'(L)=0\}$ where $D_x^2=\frac{d^2}{dx^2}.$ Then
  ${\bf{T}}$ is a densely defined, self-adjoint, positive definite, and unbounded operator on $(L^2(\Omega))^m,$  and therefore $({\bf{h}}_E {\bf T} + P {\bf{ G}}_E)^{-1}$ exists and is a bounded operator defined on all of $(L^2(\Omega))^m$ where $P={\bf{B}}{\bf{E}}_\mc{O}^{-1} {\bf{h}}_\mc{O}^{-1} {\bf{B}}^T>0.$ Now define  the operator $J=-({\bf{h}}_E {\bf T} + P {\bf{ G}}_E)^{-1}~{\bf T}.$ Then $J$ extends to a  continuous and self-adjoint  operator  on $(L^2(\Omega))^{m},$ and
  \begin{eqnarray}
  \label{eq90}J=-{\bf{h}}_E^{-1}+{\bf{h}}_E^{-1}P {\bf{ G}}_E ({\bf{h}}_E {\bf T} + P {\bf{ G}}_E)^{-1} ~~{\rm on}~~\rm{Dom}({\bf T}).
  \end{eqnarray}
  To show this, let $s=Jz=-({\bf{h}}_E {\bf T} + P {\bf{ G}}_E)^{-1}{\bf T} z$ so that $s\in {\rm Dom}({\bf T})$ and ${\bf T}z=-{\bf{h}}_E {\bf T}s + P {\bf{ G}}_E s.$ Then
  \begin{eqnarray}
 \label{eq9} ({\bf{h}}_E {\bf T} + P {\bf{ G}}_E)s=-{\bf T} z=-({\bf{h}}_E ^{-1}\left[\left({\bf{h}}_E {\bf T} + P {\bf{ G}}_E\right)z- P {\bf{ G}}_Ez\right].
     \end{eqnarray}
   By applying $({\bf{h}}_E {\bf T} + P {\bf{ G}}_E)^{-1}$ to both sides of (\ref{eq9}), we get (\ref{eq90}).\\

   Now we show that $J$ is non-positive on $(L^2(\Omega))^m.$ Let $w=({\bf{h}}_E {\bf T} + P {\bf{ G}}_E)^{-1} z$ so that $w\in {\rm Dom}({\bf T})$ and $z={\bf{h}}_E {\bf T}w + P {\bf{ G}}_E w.$ Then
 \begin{eqnarray}
 \nonumber  (Jz, z)_m &=& \int_{\Omega^m}  \left(-{\bf{h}}_E^{-1}\left({\bf{h}}_E {\bf T}w + P {\bf{ G}}_E w \right)+ {\bf{h}}_E^{-1}P {\bf{ G}}_E w \right)\left({\bf{h}}_E {\bf T}\bar w + P {\bf{ G}}_E \bar w\right)~dx\\
 \nonumber &=& -{\bf{h}}_E ({\bf T} w, {\bf T} \bar w)_m- P {\bf{ G}}_E ({\bf D}w, {\bf D}\bar w)_m\le 0
 \end{eqnarray}
where ${\bf D}={\rm{diag}}~(D_x,\ldots, D_x)$ is an operator defined on the domain
$\rm{Dom}({\bf D})=\{\psi\in (H^1(\Omega))^{m}~:~ \psi(0)=0\}$ and $D_x=\frac{d}{dx}.$ Now, we are in the position of solving (\ref{eq}). Multiplying the second equation in (\ref{eq}) by ${\bf{B}}{\bf{E}}_\mc{O}^{-1} {\bf{h}}_\mc{O}^{-1},$ and using the first and the third equations in (\ref{eq})  yields
 \begin{eqnarray}
\nonumber  &&K u'''' + N^{\rm T} {\bf{h}}_E {\bf{G}}_E {\bf{h}}_E N \left(J u'\right)'=0. 
 \end{eqnarray}
But since $J$ is non-positive, the operator $ K D_x^4 + D_x(N^{\rm T} {\bf{h}}_E {\bf{G}}_E {\bf{h}}_E N J D_x)$ is a positive operator (using (\ref{bc})). This implies that $u=0.$ Therefore $v_{\mc O}=0$ by (\ref{eq}). \qed


We have the following theorem for the well-posedness of the Cauchy problem (\ref{main})-(\ref{initial}).

\begin{theorem} \label{eigens}$\mc{A}: \mc{D}( \mc{A}) \to \mc{H}$ is the infinitesimal generator of a $C_0-$semigroup of contractions. Therefore for every $T\in \mathbb{R}^+,$ $(z^0, z_\mc{O}^0, u^1, v_\mc{O}^1)\in \mc{D}(\mc{A})$ solves (\ref{main})-(\ref{initial}), and we have
$ (z, v_\mc{O}, \dot z, \dot v_\mc{O})\in C\left([0,T]; \mc{D}(\mc{A})\right)\cap C^1\left([0,T]; \mc{H}\right).$  Moreover, the spectrum $\sigma({\mc{A}})$ of $\mc{A}$ has all isolated eigenvalues.\end{theorem}

\noindent\emph{Proof.} $\mc{A}$ is an m-dissipative operator by Lemmata \ref{surjective} and \ref{disspppp}. Therefore, $\mc{A}: \mc{D}( \mc{A}) \to \mc{H}$ is the infinitesimal generator of a $C_0-$semigroup of contraction by L\"{u}mer-Phillips theorem \cite{Pazy}.  By using the fact that $\mc{D}(\mc{A})$ is densely defined and compact in $\mc{H},$ and $0\in \rho({\mc{A}}) $ by Lemma \ref{egigen0},  $(\lambda I -\mc{A})^{-1}$ is compact at $\lambda=0,$ thus compact for all $\lambda\in \rho(\mc{A}).$ Hence the spectrum of $\mc{A}$ has all isolated eigenvalues. $\hfill\square$


\section{Uniform stabilization of the decoupled system, i.e. ${\bf{ G}}_E\equiv 0$} In this section, we prove the exponential stability of the decoupled system:
\label{Decoupled}
\begin{eqnarray}\left\{ \begin{array}{l}
\ddot z -\alpha \ddot z'' +  K   z''''  = 0 ~~~ {\rm{on}}~~ \Omega\times \mathbb{R}^+ \\
   {\bf{h}}_\mc{O} {\bf{p}}_\mc{O} {\ddot v}_\mc{O} -{\bf{h}}_\mc{O} {\bf{E}}_\mc{O}  {v}_\mc{O}''  = 0~~ {\rm{on}} ~~ \Omega\times \mathbb{R}^+\\
 \end{array} \right.
\label{maindec}
\end{eqnarray}
with  initial and boundary conditions
\begin{eqnarray}\left\{ \begin{array}{l}
\label{mainbdconds}
 z(0,t)=z'(0,t)= z(L,t)=0, z''(L,t) + \gamma_0\dot z'(L,t)=0   \\
{v}_\mc{O}(0,t) = 0,  ~ {v}_\mc{O}'(L,t) + {{\Upsilon}}_\mc{O}\dot{v}_\mc{O}(L,t)=0   \\
 z(x,0)=z^0,  ~\dot z(x,0)=z^1, ~ {v}_\mc{O}(x,0)= {v}^0_\mc{O}, ~ {\dot v}_\mc{O}(x,0)= {v}^1_\mc{O}.
  \end{array} \right.
\end{eqnarray}

\subsection{Semigroup formulation}
Let $$U=:(u,{\bf{u}})=(z, {v}_\mc{O})^{\rm T}, ~~~~V:= (v, {\bf{v}})^{\rm T}=(\dot z,  {\dot v}_\mc{O})^{\rm T}, ~~{\rm{and}} ~~ Y:=(U,V)^{\rm T}.$$  Then the semigroup corresponding to (\ref{maindec}) is given by
\begin{eqnarray}&&\frac{dY}{dt}=\mc{A}_d Y:=\left( {\begin{array}{*{20}c}
   {0} & {I }  \\
   {A_d } & {0} &   \\
\end{array}} \right) \left( \begin{array}{l}
 U \\
 V \\
 \end{array} \right), \quad Y(0)=(z^0,  {v}^0_\mc{O}, z^1,  {v}^1_\mc{O} )^{\rm T}
\label{semigroupfordec1}\end{eqnarray}
 where $A_d U:=\left( {\begin{array}{*{20}c}
   -K L^{-1} u''''  \\
   {  {\bf{p}}^{-1}_\mc{O} } {  {\bf{E}}_\mc{O} }\textbf{u}'' \\
\end{array}} \right)\label{A1dec}.$ Define the bilinear forms $a_d$ and $c_d$ by
\begin{eqnarray}
\nonumber  c_d(z, v_\mc{O}; \hat{z},  {\hat{v}_\mc{O}})&=& \left<z,\hat{z}\right>_{\Omega}+\alpha\left<z',\hat{z}'\right>_{\Omega}+ \left< {\bf{h}}_\mc{O} {\bf{p}}_\mc{O} { v}_\mc{O},   {\hat{ v}}_\mc{O}\right>_{\Omega}\\
  \nonumber a_d(z, v_\mc{O}; \hat{z},  {\hat{v}_\mc{O}})&=& K  \left<z'',\hat{z}''\right>_{\Omega}+ \left< {\bf{h}}_\mc{O} {\bf{E}}_\mc{O} { v}_\mc{O}',   { \hat{v}}_\mc{O}'\right>_{\Omega}
\end{eqnarray}
where  $a_d(\cdot), c_d(\cdot)$ are the quadratic forms that agree with $a(\cdot ; \cdot), c(\cdot ; \cdot)$ on the diagonal.

The corresponding  energy inner product on $\mc H$ is given by
\begin{eqnarray}
\nonumber \left<Y,\widehat{Y}\right>_{\mc{H}}= a_d(U;\widehat{U})+c_d(V;\widehat{V}). & 
\end{eqnarray}
\begin{theorem} \label{eigensss} The operator $\mc{A}_d: \mc{D}( \mc{A}_d) \to \mc{H}$ defined in (\ref{semigroupfordec1}) is the infinitesimal generator of a $C_0-$semigroup of contractions. Therefore for every  $(z^0,  {v}^0_\mc{O}, z^1,  {v}^1_\mc{O} )^{\rm T}\in \mc{D}(\mc{A}_d)$,  $Y$    solves (\ref{maindec})-(\ref{mainbdconds}) and
$ Y\in C\left([0,\infty); \mc{D}(\mc{A}_d)\right)~\cap~ C^1\left([0,T]; \mc{H}\right).$  Moreover, the spectrum of $\mc{A}_d$ consists of  isolated eigenvalues.\end{theorem}

\noindent\emph{Proof.} Note that $\mc{D}(\mc{A}_d)=\mc{D}(\mc{A})$ where $\mc D(\mc A)$ is defined by (\ref{eq21}).  The proof of the Theorem \ref{eigens}  remains valid when ${\bf{G}}_E\equiv 0$, and hence Theorem \ref{eigensss} follows. $\hfill\square$    

Now we find the adjoint operator $\mc{A}_d^*$ which is needed in the proof of  Lemma \ref{lemomega}.

 \begin{lemma} \label{skew-adjoint}The infinitesimal generator $\mc{A}_d$  satisfies
$$\mc{A}_d^*(\gamma)=-\mc{A}_d(-\gamma), ~~\mbox{on} ~~ \mc{D}(\mc{A}_d(-\gamma))= \mc{D}(\mc{A}_d^*)$$
where $\mc{A}_d(\gamma)$ denotes the dependence of $\mc{A}_d$ on the feedback gains $ \gamma =(\gamma_0, \gamma_1,\cdots,\break \gamma_{2m+1})$.
\end{lemma}

\noindent\emph{Proof.} Let $U_1=[u_1,{\bf{u}}_1, v_1, {\bf{v}_1}]^{\rm T}\in \mc{D}(\mc{A}_d),$ $U_2=[u_2,{\bf{u}}_2, v_2, {\bf{v}_2}]^{\rm T}\in \left(\mc{H} \cap C^{\infty}(\Omega)\right).$ Then, $U_1$ and $U_2$ satisfy the following boundary conditions
\begin{eqnarray}
\nonumber &&u_1(0)=u_1'(0)=u_1(L)=  u_2(0)=u_2'(0)=u_2(L)=0,~{\bf u}_1(0)={\bf u}_2(0)=0\\
\label{bcs} &&u_1''(L)+\gamma_0 v_1'(L)=0, ~ {\bf{u}}_1'(L)+{{\Upsilon}}_\mc{O} {\bf v}_1(L)=0.
\end{eqnarray}
A  calculation using (\ref{bcs}) shows that
\begin{eqnarray}
\nonumber  \left<\mc{A}_d(\gamma)U_1, U_2\right>_{\mc{H}} &=&  \left<U_1, -\mc{A}_d(-\gamma)U_2\right>_{\mc{H}} + K v_1'(L)\left(\gamma v_2'(L) - u_2''(L)\right)\\
\nonumber &&\quad+ {\bf{h}}_\mc{O} {\bf{E}}_\mc{O}{\bf{v}}_1(L) \cdot \left({\bf{u}}_2'(L) - {{\Upsilon}}_\mc{O} {\bf{v}}_2(L) \right).
\end{eqnarray}
This implies that $\mc{A}_d^*=-\mc{A}_d(-\gamma)$ on the space $\mc{D}(\mc{A}_d(-\gamma))$. It follows from the Lemma \ref{surjective}  that
$\mc{A}_d(-\gamma)$  has no larger closed extension and hence $\mc D(\mc A_d(-\gamma)) =\mc D(\mc A_d^*)$.
  $\hfill\square$



\subsection{Spectral analysis}
In this section, we prove  the Riesz basis property for the eigenfunctions of the   Rayleigh beam equation with
boundary feedback.  A similar analysis  applies to the wave equations with boundary feedback of the form (\ref{mainbdconds})
(see Theorem \ref{oth}).

The eigenvalue problem corresponding to Rayleigh beam in (\ref{maindec}) is given as the following
\begin{equation}\left\{ \begin{array}{l}
\label{eprob}
 K u'''' - \alpha\lambda^2  u'' +  \lambda^2u=0 \\
 u(0)=u'(0)=u(L)=0,~~ u''(L)+\gamma_0 \lambda u'(L)=0.
  \end{array} \right.
\end{equation}

Now let $\lambda=is_0.$ Then the solution  of (\ref{eprob}) is
\begin{eqnarray}
\label{gensol}& u(x)=C_1 \sin{\sqrt{\theta_0}x}+C_2 \cos{\sqrt{\theta_0 }x}+C_3 \sinh{\sqrt{\xi_0 }x}+C_4 \cosh{\sqrt{\xi_0 }x}&
\end{eqnarray}
where  \begin{eqnarray}\label{eq31}&&\theta_0(s_0) = \frac{\alpha s_0^2 +\alpha s_0^2 \sqrt{1+ \frac{4K}{\alpha^2 s_0^2}}}{2K}, ~~~ \xi_0(s_0)= \frac{\alpha s_0^2 \sqrt{1+ \frac{4K}{\alpha^2 s_0^2}}-\alpha s_0^2 }{2K}.
\end{eqnarray}
By using the first three boundary conditions $u(0)=u'(0)=0,$ and $u''(L)+is_0\gamma_0 u'(L)=0$ for  (\ref{gensol})  we get
\begin{eqnarray}
\nonumber u(x) &=&-\xi_0 \sqrt{\theta_0} \sinh{\sqrt{\xi_0}(L-x)} +\theta_0 \sqrt{\xi_0}\sin{\sqrt{\theta_0}(L-x)} \\
 \nonumber && - i\gamma s \sqrt{\xi_0\theta_0}\cosh{\sqrt{\xi_0}(L-x)}-i\gamma s_0 \sqrt{\xi_0\theta_0} \cos{\sqrt{\theta_0}(L-x)}  \\
 \nonumber &&- \theta_0\sqrt{\xi_0}\sin{\sqrt{\theta_0}L} \cosh{\sqrt{\xi_0}x} +\xi_0\sqrt{\theta_0}\sinh{\sqrt{\xi_0}L} \cos{\sqrt{\theta_0}x}  \\
 \nonumber && + \theta_0\sqrt{\theta_0}\cos{\sqrt{\theta_0}L}\sinh{\sqrt{\xi_0}x}  -\xi_0\sqrt{\xi_0}\cosh{\sqrt{\xi_0}L}\sin{\sqrt{\theta_0}x}\\
\nonumber && + is_0\gamma\sqrt{\xi_0\theta_0}  \left[ \cos{\sqrt{\theta_0}L} \cosh{\sqrt{\xi_0}x}+\cosh{\sqrt{\xi_0}L}\cos{\sqrt{\theta_0}x}\right]\\
 \label{uofx} && + is_0 \gamma\left[\theta_0 \sin{\sqrt{\theta_0}L}\sinh{\sqrt{\xi_0}x}-\xi_0\sinh{\sqrt{\xi_0}L} \sin{\sqrt{\theta_0}x}\right].
\end{eqnarray}
By using the last boundary condition $u'(L)=0,$  we obtain the characteristic equation that $s_0$ satisfies
\begin{eqnarray}
 \nonumber & -\theta_0\sqrt{\xi_0} \cosh{\sqrt{\xi_0}L}\sin{\sqrt{\theta_0}L}  -\xi_0 \sqrt{\xi_0} \cosh{\sqrt{\xi_0}L}\sin{\sqrt{\theta_0}L} &\\
\nonumber & +  \theta_0\sqrt{\theta_0} \sinh{\sqrt{\xi_0}L}\cos{\sqrt{\theta_0}L}-\xi_0\sqrt{\xi_0}\cosh{\sqrt{\xi_0}L} \sin{\sqrt{\theta_0}L}& \\
 \label{eq43} & + 2is_0\gamma\sqrt{\xi_0\theta_0}  \cos{\sqrt{\theta_0}L} \cosh{\sqrt{\xi_0}L}+ is_0  \gamma (\theta_0-\xi_0) \sin{\sqrt{\theta_0}L}\sinh{\sqrt{\xi_0}L}=0.\quad\quad.&
\end{eqnarray}
Since we have $\theta_0\xi_0 =\frac{s^2}{K},$ by (\ref{eq31}), we find that
\begin{eqnarray}\theta_0 &=& \frac{\alpha s_0^2 +\alpha s_0^2 \sqrt{1+ \frac{4K}{\alpha^2 s_0^2}}}{2K}=\frac{\alpha s_0^2 }{K}+ \frac{1}{\alpha}+O(\frac{1}{s_0^2}), ~~~~ {\mbox{as}}~~ s_0\to\infty,\label{siamtheta} \\
\label{siamxi} \xi_0&=& \frac{\alpha s_0^2 \sqrt{1+ \frac{4K}{\alpha^2 s_0^2}}-\alpha s_0^2 }{2K}= \frac{1}{\alpha}+ O(\frac{1}{s_0^2})~~~~ {\mbox{as}}~~ s_0\to\infty.\end{eqnarray}
 Multiplying (\ref{eq43}) by $\frac{1}{\theta_0\sqrt{\theta_0}\sinh{\sqrt{\xi_0}L}}$ and eventually using  (\ref{siamtheta}) and (\ref{siamxi}) yields
\begin{eqnarray}
\label{char3} \cos{\sqrt{\theta_0}L} + i \gamma s_0 \theta_0 \sinh{\sqrt{\xi_0}L} = 0, ~~~~ {\mbox{as}}~~ s_0\to\infty~~~~ ({\rm or} ~~\theta_0\to\infty).
\end{eqnarray}
Solving (\ref{char3})  is equivalent to solving
\begin{eqnarray}
\label{eq46}& e^{2i\sqrt{\theta_0}L}-\frac{\gamma_0\sqrt{\frac{K}{\alpha}}-1}{\gamma_0\sqrt{\frac{K}{\alpha}}{}+1}=O(\frac{1}{\theta_0}), ~~ {\rm as }~~ s_0\to\infty .&
\end{eqnarray}


The following theorem characterizes the eigenvalues of (\ref{eprob}).

\begin{theorem} \label{expresofevalues} Assume (\ref{assump}).  The eigenvalues $\{\lambda^{\pm}_{0,n}\}$ of (\ref{eprob})  for sufficiently large $n$ consist of
complex conjugate pairs  $\lambda^{-}_{0,n}$  and  $\lambda^{+}_{0,n}$  with  asymptotic form
$\lambda^+_{0,n}= i\sqrt{\frac{K}{\alpha}}\sigma_{0,n} + O(\frac{1}{n})$ as $n\to\infty$  where
\begin{equation}\sigma_{0,n}=\left\{ \begin{array}{l}
\nonumber 
  \frac{i}{2L}\ln{\left|\frac{\gamma_0\sqrt{\frac{K}{\alpha}}+1}{\gamma_0\sqrt{\frac{K}{\alpha}}-1}\right|}+ \frac{n\pi}{L},  \quad\quad\quad\gamma_0> \sqrt{\frac{\alpha}{K}} \\
  \frac{i}{2L}\ln{\left|\frac{\gamma_0\sqrt{\frac{K}{\alpha}}+1}{\gamma_0\sqrt{\frac{K}{\alpha}}-1}\right|}+ \frac{(n+\frac{1}{2})\pi}{L},\quad~  \gamma_0< \sqrt{\frac{\alpha}{K}}.
 \end{array} \right.
\end{equation}

\end{theorem}

\noindent\emph{Proof.}  First, note that $\{\sigma_{0,n}\}_{n\in \mathbb{Z}_+}$ are the solutions of   (\ref{eq46}) when the right hand side of the equation is zero.   We claim that $\{\sigma_{0,n} + O(\frac{1}{n})\}$  solve (\ref{eq46})   for all sufficiently large $n\in \mathbb{Z}_+.$   Without loss of generality, we only consider the case $\gamma_0> \sqrt{\frac{\alpha}{K}}.$  Let $ f(\theta_0)=e^{2iL\sqrt{\theta_0}}-\frac{\gamma_0\sqrt{\frac{K}{\alpha}}-1}{\gamma_0\sqrt{\frac{K}{\alpha}}{}+1}$ and  $g(\theta_0)=O(\frac{1}{\theta_0}).$ Now consider the circle
 $B_n=\{\theta=|\theta_0|e^{i\vartheta_0}~:~ |\sqrt{\theta_0}-\sigma_{0,n}|\le \frac{1}{|\sigma_{0,n}|^2} \}.$
 Then on the disc  $D_n=\{\theta_0~:~ |\sqrt{\theta_0}-\sigma_{0,n}|= \frac{1}{|\sigma_{0,n}|^2} \}$ we have
 \begin{eqnarray}\label{eq52}|\sqrt{\theta_0}|=|\sigma_{0,n}| \left(1+ O\left(\frac{1}{|\sigma_{0,n}|^3}\right)\right).\end{eqnarray}
By a simple calculation one easily obtains the following  for sufficiently large $n:$
 \begin{eqnarray}
 \nonumber &&\left|~ e^{2iL\sqrt{\theta_0}}-\frac{\gamma_0\sqrt{\frac{K}{\alpha}}-1}{\gamma_0\sqrt{\frac{K}{\alpha}}+1}\right|= \left|~e^{2iL\sqrt{\theta_0}}-e^{2iL\sigma_{0,n}}\right|\\
\nonumber &&=\quad \left|\frac{\gamma_0\sqrt{\frac{K}{\alpha}}-1}{\gamma_0\sqrt{\frac{K}{\alpha}}{}+1}\right|\left(\frac{2L}{|\sqrt{\theta_0}|}+O\left(\frac{1}{|\theta_0|}\right)\right) > O(\frac{1}{|\theta_0|})\quad
\end{eqnarray}
where we have used (\ref{eq43}), (\ref{eq46}) and (\ref{eq52}).      Therefore, by Rouch\'{e}'s theorem, $f(\theta_0)+g(\theta_0)$ has a unique zero in the ball $B_n$ for sufficiently large $n.$ That is, there exists a unique solution of the equation (\ref{eq46}) in $B_n.$    This proves our claim and hence $\lambda_{0,n}^+= i\sqrt{\frac{K}{\alpha}}\sigma_{0,n} + O(\frac{1}{n})$ as $n\to\infty$ by (\ref{siamtheta}). $\hfill\square$

\begin{theorem} \label{asymptots} Assume (\ref{assump}) holds.   The eigenfunctions \begin{equation*}{\left\{(e_{0,n},\lambda_{0,n} e_{0,n})^{\rm T}, ~ \forall n \in \mathbb{Z}\right\}}\end{equation*} corresponding to the eigenvalues $\{ \lambda_{0,n}\}$ of (\ref{eprob}) have the following asymptotic expressions:
\begin{eqnarray}\label{eigenfunc1}
\left( \begin{array}{c}
 {e''_{0,n}}  \\
 \lambda_{0,n} {e'_{0,n}}  \\
 \end{array} \right)=\left( \begin{array}{c}
 \cos{\sigma_{0,n} x} +O(1/n) \\
i\sqrt{ \frac{K}{\alpha}}    \sin{\sigma_{0,n} x}  +O(1/n)\\
  \end{array} \right)\end{eqnarray}
which are asymptotically  normalized in $H^2_\#(\Omega)\times H^1_0(\Omega)$.
\end{theorem}

\noindent\emph{Proof.} First, we find estimates for ${e''_{0,n}}$ and $\lambda_{0,n}{e'_{0,n}}.$
Application of  (\ref{siamtheta}),(\ref{siamxi}) and  Theorem \ref{expresofevalues} shows that for any $x\in [0,L]$ we have
\begin{eqnarray}
\nonumber \sinh{\sqrt{\xi_{0,n}}x} = \sinh{\frac{x}{\sqrt{\alpha}}} + O(\frac{1}{n}),\quad  \cosh{\sqrt{\xi_{0,n}}x} = \cosh{\frac{x}{\sqrt{\alpha}}} + O(\frac{1}{n})\\
\label{eqq21} \sin{\sqrt{\theta_{0,n}}x} = \sin{\sigma_{0,n} x} + O(\frac{1}{n}), \quad
 \cos{\sqrt{\theta_{0,n}}x} = \cos{\sigma_{0,n} x} + O(\frac{1}{n}).
\end{eqnarray}
Therefore    the characteristic equation (\ref{char3}) takes the form
\begin{eqnarray}\label{eq65}
\cos{\sqrt{\theta_{0,n}}L}+i\gamma_0\sqrt{\frac{K}{\alpha}} \sin{\sqrt{\theta_{0,n}}L}=O(\frac{1}{n}).
\end{eqnarray}
Now we use (\ref{uofx}), (\ref{siamtheta}), (\ref{siamxi}), (\ref{eqq21}) and (\ref{eq65}) to get
\begin{eqnarray}
\label{eq70}&& \frac{1}{\theta_{0,n}\sqrt{\xi_{0,n}\theta_{0,n}}}\nonumber u'_{n}(x) = \zeta_n \sin{\sqrt{\theta_{0,n}}x}  + O(\frac{1}{n}),\quad\quad
\end{eqnarray}
\begin{eqnarray}
\nonumber &\frac{1}{\theta_{0,n}^2\sqrt{\xi_{0,n}}}u''_n(x) = \zeta_n \cos{\sqrt{\theta_{0,n}}x}+ O(\frac{1}{n})&
\end{eqnarray}
where $\zeta_n=\left(\left( \frac{K\gamma_0^2}{\alpha}-1\right)\sin{\sqrt{\theta_{0,n}}L} -i\gamma_0  \sqrt{\frac{K}{\alpha}} \cosh{\sqrt{\xi_{0,n}}L}\right).$ In the following we verify that $\zeta_n=O(1)\ne 0.$ We have two cases. \vspace{0.05in}

 \noindent \underline{\textbf{Case I.}} Let $\gamma_0> \sqrt{\frac{\alpha}{K}}.$ By (\ref{siamtheta})  we have
\begin{eqnarray}
\nonumber  \zeta_n  = i(-1)^{n}\left(\frac{K}{\alpha}\gamma_0^2-1\right) \sinh{\left(\frac{1}{2}\ln{\left|\frac{\gamma_0\sqrt{\frac{K}{\alpha}}+1}{\gamma_0\sqrt{\frac{K}{\alpha}}-1}\right|}\right)} - i\gamma_0 \sqrt{\frac{K}{\alpha}}    \cosh{\sqrt{\xi_{0,n}    }L}+ O(\frac{1}{n}).\quad&
\end{eqnarray}
When  $n$ is odd, it is clear that the complex number $$\left(\left(\frac{K}{\alpha}\gamma_0^2-1\right)\sin{\left(\frac{1}{2}\ln{\left|\frac{\gamma_0\sqrt{\frac{K}{\alpha}}+1}{\gamma_0\sqrt{\frac{K}{\alpha}}-1}\right|}\right)} + i\gamma_0 \sqrt{\frac{K}{\alpha}}    \cosh{\sqrt{\xi_{0,n}    }L}\right)$$ has a nonzero (but constant) real part. Therefore $\zeta_n=O(1)\ne 0.$
\vspace{0.05in}

When $n$ is even, we must have
$$\zeta_n=i\left[\left(\frac{K}{\alpha}\gamma_0^2-1\right) \sinh{\left(\frac{1}{2}\ln{\left|\frac{\gamma_0\sqrt{\frac{K}{\alpha}}+1}{\gamma_0\sqrt{\frac{K}{\alpha}}-1}\right|}\right)} - \gamma_0 \sqrt{\frac{K}{\alpha}}    \cosh{\sqrt{\xi_{0,n}    }L}+ O(\frac{1}{n})\right].$$
But $\zeta_n=O(1)\ne 0,$ otherwise   by (\ref{siamxi}) and (\ref{eq65}), we have
\begin{eqnarray}
\nonumber   &&\gamma_0 \sqrt{\frac{K}{\alpha}}    \cosh{\frac{L}{\sqrt{\alpha}}}= \left(\frac{K}{\alpha}\gamma_0^2-1\right) \sinh{\left(\frac{1}{2}\ln{\left|\frac{\gamma_0\sqrt{\frac{K}{\alpha}}+1}{\gamma_0\sqrt{\frac{K}{\alpha}}-1}\right|}\right)}=\sqrt{\frac{K \gamma_0^2-\alpha}{{\alpha}}},
\end{eqnarray}
and by taking the square of both sides and after cancelations we get a contradiction:
\begin{eqnarray}
\nonumber &&  \frac{K \gamma_0^2}{\alpha}   \cosh^2{\frac{L}{\sqrt{\alpha}}}=\frac{K \gamma_0^2-\alpha}{{\alpha}}~~~ {\rm or} ~~~~ -\frac{\alpha}{K\gamma_0^2}= \sinh^2{\frac{L}{\sqrt{\alpha}}}.
\end{eqnarray}
\vspace{0.05in}

 \noindent \underline{\textbf{Case II.}} Now let $\gamma_0< \sqrt{\frac{\alpha}{K}}.$ Then
\begin{eqnarray}
\nonumber \zeta_n = \left(\frac{K}{\alpha}\gamma_0^2-1\right) \cosh{\left(\frac{1}{2}\ln{\left|\frac{\gamma_0\sqrt{\frac{K}{\alpha}}+1}{\gamma_0\sqrt{\frac{K}{\alpha}}-1}\right|}\right)} (-1)^{n}- i\gamma_0 \sqrt{\frac{K}{\alpha}}    \cosh{\sqrt{\xi_{0,n}    }L}+ O(\frac{1}{n}).
\end{eqnarray}
But since the last expression above  has a nonzero (but constant) real part, we have  $\zeta_n=O(1)\ne 0.$
\vspace{0.05in}

This proves the first part of the theorem.

Now if we set $$ e_{0,n}=\zeta_n^{-1}\theta_{0,n}^{-2}\xi_{0,n}^{-1/2} u_n, $$
then the second part of the theorem follows that eigenvectors are asymptotically normalized, i.e. $\|(e_{0,n}, \lambda_{0,n} e_{0,n})^{\rm T}\|_{H^2_\#(\Omega)\times H^1_0(\Omega)}\asymp 1.$ $\hfill\square$

\begin{lemma} \label{lemomega} Assume (\ref{assump}). The generalized eigenfunctions $\{(e_n,  \lambda_{0,n}e_n)^{\rm T}, ~n\in\mathbb{Z} \}$ of (\ref{eprob}) corresponding to the eigenvalues $\{\lambda_{0,n}, ~n\in \mathbb{Z}\}$ are  $\omega-$linearly independent in $ H^2_\#(\Omega)\times H^1_0(\Omega).$ Moreover, $\{(e''_n, \lambda_{0,n}e'_n, )^{\rm T}, ~n\in\mathbb{Z} \}$ is $~\omega-$linearly independent in $(L^2(\Omega))^{2}.$
\end{lemma}

\noindent\emph{Proof.} We first prove that the eigenfunctions $\{(e^*_n,  \lambda^*_{0,n}e^*_n)^{\rm T}, ~n\in\mathbb{Z} \}$ of the adjoint eigenvalue problem (see (\ref{eadjprob}) below), are biorthogonal to the eigenfunctions $\{(e_n, \lambda_{0,n}e_n)^{\rm T}, ~n\in\mathbb{Z} \}$ of (\ref{eprob}). By using Lemma \ref{skew-adjoint}, we consider the following adjoint eigenvalue problem:
\begin{equation}\left\{ \begin{array}{l}
\label{eadjprob}
 K u'''' - \alpha{\lambda^*}^2  u'' + {\lambda^*}^2 u=0\\
 u(0)=u'(0)=u(L)=0,~~ u''(L)+\gamma_0 \lambda^* u'(L)=0.
 \end{array} \right.
\end{equation}

This is exactly the same boundary value problem as (\ref{eprob}). Therefore $\lambda^*_{0,n}=\lambda_{0,n}.$ The only difference is the expression of the eigenfunctions of (\ref{eadjprob})  given by $$\{(e^*_n, \lambda^*_{0,n}e^*_n)^{\rm T}, ~n\in\mathbb{Z} \}= \{(e_n, -\lambda_{0,n}e_n, )^{\rm T}, ~n\in\mathbb{Z} \}.$$
It is possible to check that
$$\left<(e_n,  \lambda_{0,n}e_n)^{\rm T}, (e_m,  -\lambda_{0,m}e_m)^{\rm T}\right>_{H^2_\#(\Omega)\times H^1_0(\Omega)}=0$$ if $m\ne -n.$ If $m=-n,$ the inner product
\begin{eqnarray}\label{21}\left<(e_n,  \lambda_{0,n}e_n)^{\rm T}, (e_m,  -\lambda_{0,m}e_m)^{\rm T}\right>_{H^2_\#(\Omega)\times H^1_0(\Omega)}\end{eqnarray}
is uniformly bounded (from (\ref{eigenfunc1})  each term has a uniform asymptotic bound). Hence a uniform bound for (\ref{21}) exists.  Therefore, $\{(e_n, \lambda_{0,n}e_n)^{\rm T}, ~n\in\mathbb{Z} \}$ is  $\omega-$linearly independent in $H^2_\#(\Omega)\times H^1_0(\Omega).$ This proves the first part of Lemma \ref{lemomega}.

To show that $\{(e''_n,  \lambda_{0,n}e'_n)^{\rm T}, ~n\in\mathbb{Z} \}$ is $~\omega-$linearly independent in $(L^2(\Omega))^{2},$ we readjust the inner product (\ref{inner}) for $H^2_\#(\Omega)\times H^1_0(\Omega)$ with an equivalent one as the following:
 \begin{eqnarray}
\nonumber \left<\left( \begin{array}{l}
 u \\
 v \\
 \end{array} \right),\left( \begin{array}{l}
 \hat u \\
 \hat v \\
 \end{array} \right)\right>_{H^2_\#(\Omega)\times H^1_0(\Omega)}= \left<v',\hat v'\right>_{\Omega}+ \left<u'', \hat u''\right>_{\Omega}.
  \end{eqnarray}

Define the map $\mc{T}:H^2_\#(\Omega)\times H^1_0(\Omega)\to (L^2(\Omega))^{2}$ by
$T (u, v)^{\rm T} = (u'', v')^{\rm T}.$
It is clear that $\mc{T}$ is an isomorphism:  $H^2_\#(\Omega)\times H^1_0(\Omega)$ to  $(L^2(\Omega))^{2}.$  Hence $\omega-$linearly independence
 is preserved.    This proves Lemma \ref{lemomega}. \qed

\begin{theorem}\label{Rieszbasisprop} Assume (\ref{assump}). The generalized eigenfunctions $\{(e_n, \lambda_{0,n}e_n)^{\rm T}, ~n\in\mathbb{Z} \}$ of (\ref{eprob}) forms a Riesz basis in $H^2_\#(\Omega)\times H^1_0(\Omega).$
\end{theorem}

\noindent\emph{Proof.} It is known that both $\{1,\cos{\frac{n\pi x}{L}}\}_{n\in \mathbb{N}}$ and $\{\sin{\frac{n\pi x}{L}}\}_{n\in\mathbb{N}}$ are orthonormal bases in  $L^2(\Omega).$ Therefore, it is easy to see that
  $$\left\{
\left( \begin{array}{c}
  \cos{\frac{n\pi x }{L}} \\
 \sin{\frac{n\pi x }{L}} \\
 \end{array} \right),\left( \begin{array}{c}
 1 \\
0 \\
 \end{array} \right),\left( \begin{array}{c}
 \cos{\frac{n\pi x }{L}}  \\
 -\sin{\frac{n\pi x }{L}}\\
 \end{array} \right)\right\}_{n\in\mathbb{N}}$$
  is also an  orthonormal basis in $(L^2(\Omega))^2.$   Now let
  \begin{eqnarray}\nonumber S=\left( {\begin{array}{*{20}c}
   \cosh{c x} & i\sinh{c x}  \\
   \delta \sinh{c x} & i\delta\cosh{c x}  \\
\end{array}} \right)
\end{eqnarray}
where
$\nonumber \delta=
  i \sqrt{\frac{K}{\alpha}}\quad ,c=\frac{1}{2L}\ln{\left|\frac{\gamma_0\sqrt{\frac{K}{\alpha}}+1}{\gamma_0\sqrt{\frac{K}{\alpha}}-1}\right|}.
$
  $S$ is a bounded linear operator from $(L^2(\Omega))^2$ to $(L^2(\Omega))^2$  since $|S|=\delta<\infty.$ For $n\in\mathbb{N}$ we have
  \begin{eqnarray}
  \nonumber  && S
\left( \begin{array}{l}
 \cos{\frac{n\pi x}{L}} \\
\sin{\frac{n\pi x}{L}} \\
 \end{array} \right)= \left( \begin{array}{l}
\cosh{(c+\frac{in\pi}{L})x} \\
 \delta\sinh{(c+\frac{in\pi}{L})x} \\
 \end{array} \right).
  \end{eqnarray}
  Therefore $ \left\{F_n= \left( \cosh{(c_0+\frac{in\pi}{L})x}, \delta_0\sinh{(c_0+\frac{in\pi}{L})x} \right)^{\rm T}
   \right\}$
 forms a Riesz basis in  $(L^2(\Omega))^{2}.$
 Now let
$\left\{G_n=(e''_n, \lambda_{0,n}e'_n)^{\rm T}, ~n\in \mathbb{Z}\right\}$
  where  $\{(e_n, \lambda_{0,n}e_n)^{\rm T},~n\in \mathbb{Z}\}$ are the eigenvectors corresponding to the eigenvalues $\{\lambda_{0,n}.\} $  Since  $\{G_n,~ n\in\mathbb{Z}\}$ is $\omega-$ linearly independent in $L^2(\Omega)^{2}$ by Lemma \ref{lemomega}, and $\|F_{n}-G_{n}\|_{(L^2(\Omega))^{2}}=O(\frac{1}{n}),$ i.e.
  \begin{eqnarray} \label{asymdiff}
\nonumber &&  \left\|~e''_n-\cosh{\left(c_0+\frac{in\pi}{L}\right)x} \right\|_{\Omega}=O(\frac{1}{n}), \quad \\&& \left\|~\lambda_{0,n}e'_n- \delta_k\sinh{\left(c_0+\frac{in\pi}{L}\right)x}\right\|_{\Omega}=O(\frac{1}{n}).
  \end{eqnarray}
  It follows from Bari's theorem \cite{Young} that  $\{G_n, ~n\in \mathbb{Z}\}$ is a Riesz basis on $(L^2(\Omega))^{2}.$ Hence $\left\{(e_n, \lambda_{0,n}e_n)^{\rm T}, ~n\in \mathbb{Z}\right\}$ forms a Riesz basis in $H^2_\#(\Omega)\times H^1_0(\Omega).$ $\hfill\square$

   The following theorem can be obtained by the same procedure.
\begin{theorem} \label{oth}Assume (\ref{assump}).  Consider
\begin{equation}\left\{ \begin{array}{l}
\label{nothing1}
   {\bf{h}}_\mc{O} {\bf{p}}_\mc{O} {\ddot v}_\mc{O} -{\bf{h}}_\mc{O} {\bf{E}}_\mc{O}  {v}_\mc{O}''  = 0~~ {\rm{on}} ~~ \Omega\times \mathbb{R}^+ \\
{v}_\mc{O}(0,t) = 0,  ~ {v}_\mc{O}'(L,t) + {{\Upsilon}}_\mc{O}\dot{v}_\mc{O}(L,t)=0   \\
~{v}_\mc{O}(x,0)= {v}^0_\mc{O}, ~ {\dot v}_\mc{O}(x,0)= {v}^1_\mc{O}.
 \end{array} \right.
\end{equation}
Then, the eigenfunctions of (\ref{nothing1})
\begin{eqnarray}\nonumber \left\{({\bf e}_{k,n}, \lambda_{k,n} {\bf e}_{k,n})^{\rm T}, k=1,3,\ldots, 2m+1, ~n\in \mathbb{Z}\right\}
 \end{eqnarray} corresponding to the branches of eigenvalues $\{ \cup \{\lambda_{k,n}\}, ~~ k=1,3,\ldots, 2m+1, ~n\in \mathbb{Z} \}$  forms a Riesz basis in $\left(H^1_*(\Omega)\right)^{(m+1)}\times (L^2(\Omega))^{(m+1)}$ where \\
${\bf e}_{k,n}=(0,\ldots, e_{k,n},\ldots,0)^{\rm T},$
$ e_{k,n}=\theta_{k,n}^{-1} \sin{\theta_{k,n} x},$
$\lambda_{k,n}= i \sqrt{\frac{E_k}{\rho_k}}\theta_{k,n}$ for all $n\in \mathbb{Z},$ and
\begin{eqnarray}\nonumber \theta_{k,n}=\left\{ \begin{array}{l}
 \frac{i}{2L}\ln{\left|\frac{\gamma_k\sqrt{\frac{E_k}{\rho_k}}+1}{\gamma_k\sqrt{\frac{E_k}{\rho_k}}-1}\right|}+ \frac{n\pi}{L},  \gamma_k> \sqrt{\frac{E_k}{\rho_k}} \\
  \frac{i}{2L}\ln{\left|\frac{\gamma_k\sqrt{\frac{E_k}{\rho_k}}+1}{\gamma_k\sqrt{\frac{E_k}{\rho_k}}-1}\right|}+ \frac{(n+\frac{1}{2})\pi}{L}, \gamma_k< \sqrt{\frac{E_k}{\rho_k}}.
\end{array} \right.
\end{eqnarray}
\end{theorem}

    \begin{theorem}\label{exp-ray} Assume (\ref{assump}). Then the semigroup generated by $\mc{A}_d$ is exponentially stable on $\mc{H}$, i.e., $ \exists M>0$ such that
\begin{eqnarray}\mc{E}(t)\le M e^{\tilde\mu t}\mc{E}(0)\label{exp}\end{eqnarray}
where $\tilde \mu=\sup\{ \mbox{Re}\,  \lambda~|~ \lambda \in \sigma(\mc{A}_d)\}<0.$
\end{theorem}

\noindent\emph{Proof.}
 The  Riesz basis property  for the Rayleigh beam equation (Theorem \ref{Rieszbasisprop}) together with  Riesz basis property   for the system of
 wave equations (Theorem \ref{oth}) imply that the eigenfunctions $\{(e_n, {\bf e}_{k,n}, \lambda_{k,n}e_n,  \lambda_{k,n} {\bf e}_{k,n})^{\rm T}, ~n\in\mathbb{Z} \}$ of the operator $\mc{A}_d$ form a Riesz basis in $\mc{H}.$     Hence, as is well known, the growth bound for the associated semigroup is determined by spectrum of the generator.   We know from Theorems \ref{expresofevalues} and \ref{oth} that the eigenvalues
 \begin{eqnarray}
 \nonumber &&\{\lambda_{k,n}, ~k=0,1,3,\ldots, 2m+1, ~n\in \mathbb{Z}\}
 \nonumber
  \end{eqnarray}
   of $\mc A_d$  have the expressions $\lambda_{0,n}= i\sqrt{\frac{K}{\alpha}}\sigma_{0,n} + O(\frac{1}{n})$ as $n\to\infty$ and $\lambda_{k,n}= i \sqrt{\frac{E_k}{\rho_k}}\theta_{k,n}$ for all $n\in \mathbb{Z}_+.$   Furthermore since ${\mc A}_d$ is dissipative,  all eigenvalues have non-positive real parts.
 Hence,  if we show that there are no eigenvalues on the imaginary axis, then the theorem is proved.      For the wave equations, this is trivial to show, and is well-known.     For the boundary conditions we have,  for the Rayleigh beam,  the possibility of imaginary eigenvalues lead to the following overdetermined eigensystem
 \begin{eqnarray}\nonumber \left\{ \begin{array}{l}
K u'''' - \alpha\lambda^2  u'' + \lambda^2u=0 \\
u(0)=u'(0)=u(L)=u'(L)= u''(L)=0.
 \end{array} \right.
\end{eqnarray}


This system was shown in \cite{O-Hansen1}  to have only the trivial solution.    Therefore $\{e^{\mc{A}_d t }\}_{t\ge 0}$ is an exponentially stable semigroup on  $\mc H,$ and (\ref{exp}) holds. $\hfill\square$

%
%

\section{Uniform stabilization of the coupled system}
\label{Coupled-Stab}
In this section, we show that one boundary feedback for each equation is enough to obtain the uniform stabilization of the multilayer RN beam. First, we will consider the decomposition $\mc{A}=\mc{A}_d+ \mc{B}$ of the semigroup generator of the original problem (\ref{semigroupfor1}) where $\mc{A}_d$ is the semigroup generator of the decoupled system and it is defined by (\ref{semigroupfordec1}), and the operator $\mc{B}:\mc{H}\to\mc{H}$ is the coupling between the layers defined as the following

\begin{eqnarray}\label{opB}\mc{B}\left( \begin{array}{c}
 u \\
 \bf u \\
 v\\
 \bf v
 \end{array} \right)=
 \left( \begin{array}{c}
 0 \\
 0_{O}\\
 L^{-1}\left(N^{\rm T} {\bf{h}}_E {\bf{G}}_E ~\phi_E'\right)   \\
-{ {\bf{h}}^{-1}_{\mc O} {\bf{p}}^{-1}_{\mc O} }{\bf{B}}^{\rm T} {\bf{G}}_E~\phi_E\\
   \end{array} \right)
 \end{eqnarray}
 where $\phi_E={\bf{h}}_E^{-1} {\bf{B}} {\bf{u}}+  N u'$

\vspace{0.1in}

\begin{lemma} \label{compact} The operator $B:\mc{H}\to \mc{H}$ defined in (\ref{opB}) is compact.
\end{lemma}

 When $(u,{\bf u}, v, {\bf v})^{\rm T}\in \mc H,$ we have $u\in H^2_\#(\Omega)$ and ${\bf u}\in (H^1_*(\Omega))^{(m+1)},$ and therefore $\phi_E\in (H^1_*(\Omega))^{(m+1)}.$ Since $L: H^2(\Omega)\cap H^1_0(\Omega)\to L^2(\Omega)$ is an isomorphism, the last terms in (\ref{opB}) satisfy
\begin{eqnarray}
\nonumber &&L^{-1}\left(N^{\rm T} {\bf{h}}_E {\bf{G}}_E ~\phi_E'\right) \in H^2(\Omega)\cap H^1_0(\Omega)\\
\label{dumb1} &&-{ {\bf{h}}^{-1}_{\mc O} {\bf{p}}^{-1}_{\mc O} }{\bf{B}}^{\rm T} {\bf{G}}_E~\phi_E \in (H^1_*(\Omega))^{(m+1)},
\end{eqnarray}
which are compactly embeddded  in $H^1_0(\Omega)$ and $(L^2(\Omega))^{(m+1)},$ respectively. Hence the operator $B$  is compact in $\mc H$. $\hfill\square$

\begin{theorem} \label{stronglystable} Assume (\ref{assump}). Then the semigroup generated by $\mc{A}$ is strongly stable in $\mc{H}.$
\end{theorem}

\noindent\emph{Proof.}
We know that our system (\ref{main})-(\ref{initial}) is dissipative by (\ref{Dissp}). If we can show that there are no eigenvalues on the imaginary axis, or in other words, the set
\begin{eqnarray}\{Y\in \mc{H} ~|~ {\rm Re} \left<\mc{A}Y, Y\right>_{\mc{H}} = -K \gamma_0 |v'(L)|^2 - {\bf{h}}_\mc{O} {\bf{E}}_\mc{O} {{\Upsilon}}_\mc{O}{\bf{v}}(L) \cdot{\bf{{\bar v}}} (L)= 0\}
\label{Eq81}
\end{eqnarray}
has only the trivial solution, i.e. $u=0, {\bf u}=0,$ then by La Salle's invariance principle, the system (\ref{main})-(\ref{initial}) is strongly stable.
Since  we have eliminated the possibility of a zero eigenvalue in  Lemma   \ref{egigen0},   (\ref{Eq81})  corresponds to
$v'(L)=\lambda u'(L)=0$ and
${\bf{v}}(L)=\lambda {\bf{u}}(L)=0$ where $\lambda\ne 0.$ Therefore, proving the strong stability of the (\ref{main})-(\ref{initial}) reduces to showing that the following eigenvalue problem
\begin{equation}\left\{ \begin{array}{l}
\lambda^2 u -\alpha \lambda^2 u'' +  K   u'''' -   N^T {\bf{h}}_E {\bf{ G}}_E \phi_E' = 0 ~~~ {\rm{on}}~~ \Omega\\
   {\bf{h}}_\mc{O} {\bf{p}}_\mc{O}\lambda^2 {\bf u}+ {\bf{h}}_\mc{O} {\bf{E}}_\mc{O}  {\bf u}'' + {\bf{B}}^T  {\bf{ G}}_E \phi_E  = 0~~ {\rm{on}} ~~ \Omega \\
 {\rm {where}}\quad ( {\bf{B}} {\bf u}={\bf{h}}_E \phi_E-{\bf{h}}_E  N u')
 \end{array} \right.
\nonumber
\end{equation}
with initial and overdetermined boundary conditions
\begin{equation}
\nonumber \left\{ \begin{array}{l}
 u(0)=u'(0)= u(L)= u'(L)=u''(L)=0    \\
{\bf u}(0) = {\bf u}(L)= {\bf u}'(L)=0
 \end{array} \right.
\end{equation}
has only the trivial solution, i.e. $u=0, {\bf u}=0.$    This same overdetermined system  came up in proving  observability for the corresponding boundary control problem
in \cite{O-Hansen3}, where the uniqueness of the zero solution was proved using a multiplier type argument.       $\hfill\square$


Now we prove our main theorem for the exponential stability of the solutions (\ref{main})-(\ref{initial}):

\noindent\emph{Proof of Theorem \ref{finaltheorem}.} We know that $\mc{A}=\mc{A}_d+\mc{B}.$ The semigroup $\{e^{(\mc{A}_d+\mc{B})t}\}_{t\ge 0}$ is strongly stable on $\mc{H}$ by Theorem \ref{stronglystable}  and the operator $\mc{B}$ is a compact in $\mc{H}$ by Lemma \ref{compact}. Therefore, since the semigroup generated by $(\mc{A}_d+\mc{B})-\mc{B}$ is uniformly exponentially stable in $\mc{H}$ then the semigroup $\{e^{(\mc{A}_d+\mc{B})t}\}_{t\ge 0}=\{e^{\mc{A}t}\}_{t\ge 0}$ is uniformly exponentially stable in $\mc{H}$ by e.g.,
the perturbation theorem of Triggiani  \cite{Trigg}. $\hfill\square$


\medskip
Received October 2012; revised April 2013.
\medskip


\begin{thebibliography}{99}


\bibitem{AAA} (MR2713103)
\newblock A. A. Allen,
\newblock \emph{Stability Results for Damped Multilayer Composite Beams and Plates},
\newblock Ph.D. thesis, Iowa State University, 2009.\vspace*{2pt}

\bibitem{AH1} (MR2679641) [10.3934/dcdsb.2010.14.1279]
\newblock A. A. Allen and S. W. Hansen,
\newblock \emph{\emph{Analyticity and optimal damping for a multilayer Mead-Markus sandwich beam,}}
\newblock \emph{Discrete Contin. Dyn. Syst. Ser. B (4),} \textbf{14} (2010), 1279--1292.\vspace*{2pt}

\bibitem{AH2} (MR2671960) [10.1016/j.na.2009.02.063]
\newblock A. A. Allen and S. W. Hansen,
\newblock \emph{\emph{Analyticity of a multilayer Mead-Markus plate,}}
\newblock \emph{Nonlinear Analysis (12),} \textbf{71} (2009), e1835--e1842.\vspace*{2pt}

\bibitem{Chen} (MR603083) [10.1137/0319008]
\newblock G. Chen,
\newblock \emph{\emph{A note on the boundary stabilization of the wave equation,}}
\newblock \emph{SIAM J. Control Optim.,}  \textbf{19} (1981), 106--113.\vspace*{2pt}

\bibitem{Fabiano-Hansen} (MR1988333)
\newblock R. H. Fabiano and S. W. Hansen,
\newblock \emph{\emph{Modeling and analysis of a three-layer damped sandwich beam,}}
\newblock \emph{Discrete Contin. Dyn. Syst.,}  (2001), Added Volume, 143--155.\vspace*{2pt}

\bibitem{Guo} (MR1892235) [10.1023/A:1017912031840]
\newblock B. Z. Guo,
\newblock \emph{\emph{Basis property of a Rayleigh beam with boundary stabilization,}}
\newblock \emph{J. Optim. Theory Appl.,}  \textbf{112} (2002), 529--547.\vspace*{2pt}

\bibitem{Hansen} (MR2078434) [10.1142/S0218202504003568]
\newblock S. W. Hansen,
\newblock \emph{\emph{Several related models for multilayer sandwich plates,}}
\newblock \emph{Math. Models Methods Appl. Sci.,}  \textbf{14} (2004), 1103--1132.\vspace*{2pt}

\bibitem{H-L} (MR1758723) [10.1142/S0218202500000306]
\newblock S. W. Hansen and I. Lasiecka,
\newblock \emph{\emph{Analyticity, hyperbolicity and uniform stability of semigroups arising in models of composite beams,}}
\newblock \emph{Math. Models Meth. Appl. Sci.,}  \textbf{10} (2000), 555--580.\vspace*{2pt}

\bibitem{Hansen-Oleg} (MR2833259) [10.3934/mcrf.2011.1.189]
\newblock S. W. Hansen and O. Y. Imanuvilov,
\newblock \emph{\emph{Exact controllability of a multilayer Rao-Nakra plate with free boundary conditions,}}
\newblock \emph{Math. Control Relat. Fields,}  \textbf{1} (2011), 189--230.\vspace*{2pt}

\bibitem{Hansen-Oleg1} (MR2859867) [10.1051/cocv/2010040]
\newblock S. W. Hansen and O. Y. Imanuvilov,
\newblock \emph{\emph{Exact controllability of a multilayer Rao-Nakra plate with clamped boundary conditions,}}
\newblock \emph{ESAIM Control Optim. Calc. Var.,}  \textbf{17} (2011), 1101--1132.\vspace*{2pt}

\bibitem{Rajaram-Hansen1} [10.1109/CDC.2005.1582645]
\newblock S. W. Hansen and R. Rajaram,
\newblock \emph{\emph{Simultaneous boundary control of a Rao-Nakra sandwich beam,}}
\newblock \emph{Proc. 44th IEEE Conference on Decision and Control and the European Control Conference,} (2005), 3146--3151.\vspace*{2pt}

\bibitem{Rajaram-Hansen3} (MR2192693)
\newblock S. W. Hansen and R. Rajaram,
\newblock \emph{\emph{Riesz basis property and related results for a Rao-Nakra sandwich beam,}}
\newblock \emph{Discrete Contin. Dyn. Syst.,} (2005), suppl., 365--375.\vspace*{2pt}

\bibitem{H-S} [10.1006/jsvi.1996.0913]
\newblock S. W. Hansen and R. D. Spies,
\newblock \emph{\emph{Structural damping in laminated beams due to interfacial slip,}}
\newblock \emph{Journal of Sound and Vibration,}  \textbf{204} (1997), 183--202.\vspace*{2pt}

\bibitem{K-Z} (MR1054123)
\newblock V. Komornik and E. Zuazua,
\newblock \emph{\emph{A direct method for the boundary stabilization of the wave equation,}}
\newblock \emph{J. Math. Pures Appl.,}  \textbf{69} (1990), 33--54.\vspace*{2pt}

\bibitem{Lagnese} (MR719445) [10.1016/0022-0396(83)90073-6]
\newblock J. Lagnese,
\newblock \emph{\emph{Decay of solutions of wave equations in a bounded region with boundary dissipation,}}
\newblock \emph{J. Differential Equations,} \textbf{50} (1983), 163--182.\vspace*{2pt}

\bibitem{LT1} (MR1122306) [10.1016/0022-0396(91)90022-2]
\newblock I. Lasiecka and R. Triggiani,
\newblock \emph{\emph{Exact controllability and uniform stabilization of Kirchhoff plates with boundary controls only in $\left. {\Delta w} \right|_\Sigma,$}}
\newblock \emph{J. Differential Equations,} \textbf{93} (1991), 62--101.\vspace*{2pt}

\bibitem{LT2} (MR876804) [10.1016/0022-0396(87)90025-8]
\newblock I. Lasiecka and R. Triggiani,
\newblock \emph{\emph{Uniform exponential decay of wave equations in a bounded region with $L_2(0,\infty; L_2(\Gamma))$-feedback control in the Dirichlet boundary conditions,}}
\newblock \emph{J. Differential Equations,} \textbf{66} (1987), 340--390.\vspace*{2pt}

\bibitem{Mead-Marcus} [10.1016/0022-460X(69)90193-X]
\newblock D. J. Mead and S. Markus
\newblock \emph{\emph{The forced vibration of a three-layer, damped sandwich beam with arbitrary boundary conditions,}}
\newblock \emph{J. Sound Vibr.,}  \textbf{10} (1969), 163--175.\vspace*{2pt}

\bibitem{O-Hansen1} (MR2852398) [10.1007/s00498-011-0069-4]
\newblock A. \"{O}. \"{O}zer and S. W. Hansen,
\newblock \emph{\emph{Exact controllability of a Rayleigh beam with a single boundary control,}}
\newblock \emph{Math. Control Signals Systems,} \textbf{23} (2011), 199--222.\vspace*{2pt}

\bibitem{O-Hansen3}
\newblock A. \"{O}. \"{O}zer and S. W. Hansen,
\newblock \emph{Exact boundary controllability results for a multilayer Rao-Nakra sandwich beam,}
\newblock to appear in SIAM J. Cont. Optim.. \vspace*{2pt}

\bibitem{Pazy} (MR710486) [10.1007/978-1-4612-5561-1]
\newblock A. Pazy,
\newblock \emph{Semigroups of Linear Operators and Applications to Partial Differential Equations,}
\newblock Applied Mathematical Sciences, 44. Springer-Verlag, New York, 1983.\vspace*{2pt}

\bibitem{Rajaram-Hansen2} (MR2332008) [10.1016/j.sysconle.2007.03.007]
\newblock R. Rajaram,
\newblock \emph{\emph{Exact boundary controllability result for a Rao-Nakra sandwich beam,}}
\newblock \emph{Systems Control Lett.,}  \textbf{56} (2007), 558--567.\vspace*{2pt}

\bibitem{B. Rao} (MR1375957) [10.1007/BF01204704]
\newblock B. Rao,
\newblock \emph{\emph{A compact perturbation method for the boundary stabilization of the Ragleigh beam equation,}}
\newblock \emph{Appl. Math. Optim.,}  \textbf{33} (1996), 253--264.\vspace*{2pt}

\bibitem{Rao-Nakra}
\newblock Y. V. K. S Rao and B. C. Nakra,
\newblock \emph{\emph{Vibrations of unsymmetrical sandwich beams and plates with viscoelastic cores,}}
\newblock \emph{J. Sound Vibr.,}  \textbf{34} (1974), 309--326.\vspace*{2pt}

\bibitem{Trigg0} (MR0445388) [10.1016/0022-247X(75)90067-0]
\newblock R. Triggiani,
\newblock \emph{\emph{On the stabilizability problem in Banach space,}}
\newblock \emph{J. Math. Anal. Appl.,}  \textbf{52} (1975), 383--403.\vspace*{2pt}

\bibitem{Trigg} (MR953013) [10.1090/S0002-9939-1989-0953013-0]
\newblock R. Triggiani,
\newblock \emph{\emph{Lack of uniform stabilization for noncontractive semigroups under compact perturbation,}}
\newblock \emph{Proc. Amer. Math. Soc.,}  \textbf{105} (1989), 375--383.\vspace*{2pt}

\bibitem{Wang-Guo} (MR2399703) [10.1007/s10957-007-9341-7]
\newblock J. M. Wang and B. Z. Guo,
\newblock \emph{\emph{Analyticity and dynamic behavior of a damped three-layer sandwich beam,}}
\newblock \emph{J. Optim. Theory Appl.,}  \textbf{137} (2008), 675--689.\vspace*{2pt}

\bibitem{Wang-Xu-Yung} (MR2193496) [10.1137/040610003]
\newblock J. M. Wang, G. Q. Xu and S. P. Yung,
\newblock \emph{\emph{Exponential stabilization of laminated beams with structural damping and boundary feedback controls,}}
\newblock \emph{SIAM J. Cont. Optim.,}  \textbf{44} (2005), 1575--1597.\vspace*{2pt}

\bibitem{Guo-Chentouf} (MR2192066) [10.1051/cocv:2005030]
\newblock J. M. Wang, B. Z. Guo and B. Chentouf,
\newblock \emph{\emph{Boundary feedback stabilization of a three-layer sandwich beam: Riesz basis approach,}}
\newblock \emph{ESAIM Control Optim. Calc. Var.,}  \textbf{12} (2006), 12--34.\vspace*{2pt}

\bibitem{Yan-Dowell}
\newblock M. J. Yan and E. H. Dowell,
\newblock \emph{\emph{Governing equations for vibratory constrained-layer damping sandwich plates and beams,}}
\newblock \emph{J. Appl. Mech.,}  \textbf{39} (1972), 1041--1046.\vspace*{2pt}

\bibitem{Young} (MR1836633)
\newblock R. Young,
\newblock \emph{An Introduction to Nonharmonic Fourier Series,}
\newblock  Revised first edition. Academic Press, Inc., San Diego, CA, 2001.


\end{thebibliography}
\end{document}